\documentclass[12pt, final, 3p, times]{Qarticle}

\usepackage{epsfig}
\usepackage{amssymb}
\usepackage{float}
\usepackage{amsfonts}
\usepackage{mathrsfs}
\usepackage{amsthm}
\usepackage{amsmath}
\usepackage{amscd}
\usepackage{lineno}

\usepackage[colorlinks,citecolor=blue,linkcolor=blue,urlcolor=blue]{hyperref}
\usepackage{graphicx}
\usepackage{subfigure}

\usepackage{stmaryrd}
\usepackage{geometry}
\usepackage{fancyhdr}
\theoremstyle{plain}

\theoremstyle{definition}

\numberwithin{equation}{section}
\selectfont
\allowdisplaybreaks[4]

\begin{document}

\title{Spectral asymptotics for linear elasticity with mixed  \   \  \  \ boundary conditions}
\fancyhead[C]{\footnotesize\textit{G. Q. Liu / Spectral asymptotics for linear elasticity with mixed boundary conditions}}

\begin{frontmatter}

\author[1]{Genqian Liu\corref{c}}
\ead{liuswgq@163.com}
 \cortext[c]{Corresponding author.}
 \address[1]{School of Mathematics and Statistics, Beijing Institute of Technology, Beijing 100081,  P. R. China}

\date{}

\begin{abstract} 

In this note, we shall show that the two-term spectral asymptotics  for the operator of linear elasticity with mixed
boundary conditions which were given by Capoferri and Mann  in \cite{CaMa-24}   essentially are  old well-known  results due to  
T. Branson, P.  Gilkey,  B. {\O}rsted and A. Pierzchalski in \cite{BGOP}.  In addition, we further point out that the so-called ``general formulae'' in \cite{SaVa-97, CaFrLeVa-23}   and the calculations for several examples in \cite{CaMa-24, CaFrLeVa-23} are all wrong.

\end{abstract}


\end{frontmatter}

\section{Introduction}

Let  $(M, g )$ be a compact connected smooth Riemannian
manifold of dimension $d \ge  2$ with smooth boundary $\partial M$.   The operator of linear elasticity $\mathscr{L}$ acting on vector fields $\mathbf{u}$ on $M$ is defined  as
 \begin{eqnarray}  \label{24.1.20-1} (\mathscr{L} \mathbf{u})^\alpha : = -\mu \big( \nabla_\beta\nabla^\beta u^\alpha +\mbox{Ric}^\alpha_{\;\;\,\beta} u^\beta \big) -(\lambda+\mu) \nabla^\alpha \nabla_\beta u^\beta, \;\;\;\, \; \alpha=1, \cdots, d. \end{eqnarray}
Here and further on we adopt the Einstein summation convention over repeated indices. The quantities $\lambda$  and $\mu$ are real constants known as Lam\'{e} parameters, assumed to satisfy the conditions
\begin{eqnarray} \label{24.1.20-2}  \mu>0, \;\;\, \lambda+2\mu>0,\end{eqnarray} 
which guarantee strong ellipticity of the operator $\mathscr{L}$ and thus  the spectral problems can be well discussed (see,  the detailed explanations in Section 2 of \cite{Liu-2309}).

Consider the following four (elastic) eigenvalue problems:
\begin{eqnarray}  \label{24.1.20-5}  \left\{ \begin{array} {ll} \mathscr{L} \mathbf{u} =\Lambda  \mathbf{u}  \;\;\; &\mbox{in} \;\;\; M,\\
 \mathbf{u}=0 \;\;\; & \mbox{on}\;\;\, \partial M, \end{array} \right.\end{eqnarray}
\begin{eqnarray}  \label{24.1.20-6}  \left\{ \begin{array} {ll} \mathscr{L} \mathbf{u} =\Lambda  \mathbf{u}  \;\;\; &\mbox{in} \;\;\; M,\\
 \mathscr{T}\mathbf{u}=0 \;\;\; &\mbox{on}\;\;\, \partial M, \end{array} \right.\end{eqnarray}
  \begin{eqnarray}  \label{24.1.20-7}  \left\{ \begin{array} {ll} \mathscr{L} \mathbf{u} =\Lambda  \mathbf{u}  \;\;\; &\mbox{in} \;\;\; M,\\  \mathbf{u} -\Big(g_{\alpha\beta} n^\alpha u^\beta\Big) \mathbf{n}  =0\;\;\;\; \,\mbox{and}\;\;\;  g_{\alpha\beta} n^\alpha \big(\mathscr{T} \mathbf{u}\big)^\beta=0\;  \;\;\; &\mbox{on}\;\;\, \partial M, \end{array} \right.\end{eqnarray}
 \begin{eqnarray}  \label{24.1.20-8}  \left\{ \begin{array} {ll} \mathscr{L} \mathbf{u} =\Lambda  \mathbf{u}  \;\;\; &\mbox{in} \;\;\; M,\\  \mathscr{T}\mathbf{u} -\Big(g_{\alpha\beta} n^\alpha (\mathscr{T} \mathbf{u})^\beta\Big) \mathbf{n}  =0\;\;\;\; \,\mbox{and}\;\;\;  g_{\alpha\beta} n^\alpha {u}^\beta=0\;  \;\;\; &\mbox{on}\;\;\, \partial M, \end{array} \right.\end{eqnarray}
where  $\mathscr{T}$ is the boundary traction operator defined as 
\begin{eqnarray} \label{24.1.20-9}  \big(\mathscr{T} \mathbf{u}\big)^\alpha := \lambda n^\alpha \nabla_\beta u^\beta +\mu \Big( n^\beta \nabla_\beta u^\alpha + n_\beta \nabla^\alpha u^\beta \Big), \;\;\; \, \, \; \alpha=1,\cdots, d. \end{eqnarray}
Here $\mathbf{n}=(n^1, \cdots, n^d)$  is the exterior unit normal vector to the boundary $\partial M$. The boundary conditions in (\ref{24.1.20-5}), (\ref{24.1.20-6}), (\ref{24.1.20-7}), (\ref{24.1.20-8})  are known as  the Dirichlet boundary condition, free boundary condition, Dirichlet-free (DF) boundary conditions and free-Dirichlet (FD) boundary conditions, respectively.  The boundary conditions DF and FD (in (\ref{24.1.20-7}) and (\ref{24.1.20-8})) are also called  {\it  mixed boundary conditions} (see,  \cite{CaMa-24}). 

In each of the four cases (\ref{24.1.20-5})--(\ref{24.1.20-8}), the spectrum is discrete and we arrange the 
eigenvalues in non-decreasing order (repeated according to multiplicity)
\begin{eqnarray*}  && 0< \Lambda_1^{\mbox{Dir} }\le  \Lambda_2^{\mbox{Dir}} \le \cdots \le \Lambda_k^{\mbox{Dir}}\le \cdots; \\
&& 0 \le  \Lambda_1^{\mbox{free} }\le  \Lambda_2^{\mbox{free}} \le \cdots \le \Lambda_k^{\mbox{free}}\le \cdots; \\
&& 0 \le  \Lambda_1^{\mbox{DF} }\le  \Lambda_2^{\mbox{DF}} \le \cdots \le \Lambda_k^{\mbox{DF}}\le \cdots; \\
&& 0 \le  \Lambda_1^{\mbox{FD} }\le  \Lambda_2^{\mbox{FD}} \le \cdots \le \Lambda_k^{\mbox{FD}}\le \cdots. \end{eqnarray*}
For each set of boundary conditions  $ \aleph\in \{\mbox{Dir}, \mbox{free}, \mbox{DF}, \mbox{FD}\}$,  the corresponding {\it eigenvalue counting function} $N_{\aleph}:\mathbb{R}\to \mathbb{N}$  is defined as 
\begin{eqnarray} \label{24.1.20-10} N_{\aleph}(\Lambda): =\#\big\{k\,\big| \,\Lambda_{k}^{\aleph}<\Lambda\big\}.\end{eqnarray}
Clearly, the function (\ref{24.1.20-10}) is monotonically non-decreasing in  $\Lambda$ and vanishes identically for $\Lambda\le \Lambda_1^{\aleph}$.

Suppose that $(\Omega, g)$ is such that the corresponding billiards is neither dead-end nor absolutely periodic.  Then there exists  two-term asymptotic expansion for the counting function $N_\aleph (\Lambda)$ : \begin{eqnarray} \label{2023.5.5-1}  {N}_\aleph (\Lambda) = a_0\, \mbox{Vol}_d (M)\, \Lambda^{d/2} + a_1^{\aleph}\, \mbox{Vol}_{d-1} (\partial M) \,\Lambda^{(d-1)/2}+ o(\Lambda^{(d-1)/2}) \;\;\, \mbox{as}\;\; \Lambda\to +\infty\;\; \end{eqnarray}
for any set of boundary conditions $\aleph \in  \{\mbox{Dir}, \mbox{free}, \mbox{DF}, \mbox{FD}\}$, 
  where  \begin{eqnarray}\label{2023.5.5-2}\!\!\!\!\!\!\!\! \!\!\!\!&& a_0= \frac{1}{(4\pi)^{d/2}\, \Gamma (1+\frac{d}{2})} \left(\frac{d-1}{\mu^{d/2}}+ \frac{1}{(\lambda+2\mu)^{d/2}}\right), \end{eqnarray}  
    ${\mbox{Vol}}_{d}(M)$ denotes the $d$-dimensional volume of $M$ and   ${\mbox{Vol}}_{d-1} (\partial M)$ denotes the $(d-1)$-dimensional volume of $\partial M$.
  Here  the leading term coefficient $a_0$  can be easily  found by calculating the kernel for the Euclidean  case  for the second-order term of $\mathscr{L}\,$ (cf. \cite{Ho4}).
  
In \cite{CaFrLeVa-23},   M. Capoferri,  L. Friedlander,  M. Levitin and D. Vassiliev gave  the second term coefficient in asymptotic formulae  (\ref{2023.5.5-1}) for the Dirichlet and free boundary conditions (see also \S6.3 of \cite{SaVa-97} for two- and three-dimensional Euclidean space cases):
\begin{eqnarray}\label{2023.5.5-3} 
\!\!\!\!\!\!\!\!\!\!\!\!\!\!\!\!&&\!\!\!\!\!\!\!\!\!\!\!\!\!\label{2024.5.5-10} a_1^{\mbox{Dir}} =- \frac{ \mu^{\frac{1-d}{2}}}{ 2^{d+1} \pi^{\frac{d-1}{2}} \Gamma(\frac{d+1}{2})} \bigg(\alpha^{\frac{d-1}{2}} +d-1 \\
 &&\quad \quad  \quad \quad  + \frac{4(d-1)}{\pi} \int_{\sqrt{\alpha}}^1 \tau^{d-2} \arctan \Big(\sqrt{(1-\alpha \tau^{-2})(\tau^{-2}-1)} \Big)\,  d\tau  \bigg), \nonumber\\
\!\!\!\!\! \!\!\!\!\!\!\!\!\!\!\!\!\!&&\!\!\!\! \!\!\!\!\!\!\!\!  \label{2024.5.5-3}  a_1^{\mbox{free}} = \frac{ \mu^{\frac{1-d}{2}}}{ 2^{d+1} \pi^{\frac{d-1}{2}} \Gamma(\frac{d+1}{2})} \bigg(\alpha^{\frac{d-1}{2}} +d-5+4 \gamma_R^{1-d}\\
&& \quad\;\;\quad\quad   \;\, + \frac{4(d-1)}{\pi} \!\int_{\sqrt{\alpha}}^1 \!\tau^{d-2} \arctan \!\Big(\frac{(\tau^{-2} \!-2)^2}{4\sqrt{(1-\!\alpha \tau^{-2})(\tau^{-2}\!-1)}} \Big) \,d\tau  \bigg),\nonumber\end{eqnarray}
 where $\alpha:= \frac{\mu}{\lambda +2\mu}$, $\gamma_R:= \sqrt{\omega_1}$, and $w_1$ is the distinguished real root of the cube equation $R_\alpha (w):=w^3-8w^2 +8(3-2\alpha) w +16 (\alpha-1)=0$ in the interval $(0,1)$.

However,  in \cite{Liu-2309}  we have given a  counter-example to show that the second coefficient  in (\ref{2024.5.5-10})--(\ref{2024.5.5-3}) are wrong. These incorrect coefficients in (\ref{2024.5.5-10})--(\ref{2024.5.5-3}) stem from  a wrong ``algorithm''  theory  which was suggested  in \cite{SaVa-97} (see also \cite{Va-84} and \cite{Va-86}). 
For the Dirichlet and free boundary condition, the correct second coefficients of the corresponding counting functions are obtained  in   \cite{Liu-21} (see also \cite{Liu-23} or \cite{Liu-22}):  
\begin{eqnarray}\label{2024.1.22-1} 
\!\!\!\!\!\!\!\!\!\!\!\!\!\!\!\!&&\!\!\!\!\!\!\!\!\!\!\!\!\!\label{202.5.5-10} a_1^{\mbox{Dir}} \,=- \frac{ \mu^{\frac{1-d}{2}}}{ 2^{d+1} \pi^{\frac{d-1}{2}} \Gamma(\frac{d+1}{2})} \left(
  \alpha^{\frac{d-1}{2}} +d-1\right),  \\
  & &\, = -\frac{1}{ 2^{d+1} \pi^\frac{d-1}{2} \Gamma\big(\frac{d+1}{2}\big)}\bigg( \frac{d-1}{\mu^\frac{d-1}{2}} + \frac{1}{(\lambda+2\mu)^\frac{d-1}{2}}\bigg)\nonumber \\  
\!\!\!\!\! \!\!\!\!\!\!\!\!\!\!\!\!\!&&\!\!\!\! \!\!\!\!\!\!\!\!  \label{2024.1.22-2}  a_1^{\mbox{free}} = \frac{ \mu^{\frac{1-d}{2}}}{ 2^{d+1} \pi^{\frac{d-1}{2}} \Gamma(\frac{d+1}{2})} \left( \alpha^{\frac{d-1}{2}} +d-1\right)\\
  & & \,\;= \frac{1}{ 2^{d+1} \pi^\frac{d-1}{2} \Gamma\big(\frac{d+1}{2}\big)}\bigg( \frac{d-1}{\mu^\frac{d-1}{2}} + \frac{1}{(\lambda+2\mu)^\frac{d-1}{2}}\bigg)\nonumber
. \end{eqnarray}
   Note that  there are not  integral terms in the correct coefficient expressions (\ref{202.5.5-10})--(\ref{2024.1.22-2}). 
   
   In \cite{CaMa-24},  M. Capoferri and I. Mann claimed that they gave the second term coefficients in  (\ref{2023.5.5-1})  for the mixed  DF and FD boundary conditions:  
 \begin{eqnarray}  \label{24.1.21-13}   \!\!\!\!\! a_1^\aleph  := @ \frac{1}{2^{d+1}  \pi^{\frac{d-1}{2}} \Gamma\big(\frac{d+1}{2}\big)} \left(\frac{d-3}{\mu^{\frac{d-1}{2}}} +\frac{1}{(\lambda+2\mu)^{\frac{d-1}{2}}}\right) \;\;  \mbox{with}\;\;   @ =\left\{\! \begin{array}{ll} - \;\;\;\mbox{for}\;\,\, \aleph =DF,\\
+ \;\;\;\; \mbox{for}\;\;\, \aleph=FD.\end{array}\right.\end{eqnarray}
There are not integral terms here either. 

In this note, we shall show that the second coefficient expressions (\ref{24.1.21-13})  of \cite{CaMa-24} for mixed boundary conditions    are old results which were given by  T. Branson, P.  Gilkey,  B. {\O}rsted and A. Pierzchalski in \cite{BGOP}. In addition, we point out a series of serious mistakes in \cite{CaMa-24} for calculating the elastic eigenvalues in the Euclidean disk and  flat cylinders. Such kinds of  fundamental  calculating errors occurred many times in the earlier  wrong  papers  \cite{CaFrLeVa-23} and  \cite{LeMoSe-21}.

 \vskip 0.82 true cm

\section{Asymptotic formulae for the mixed boundary conditions}

\vskip 0.32 true cm

Let  $(M, g )$ be a compact connected smooth Riemannian
manifold of dimension $d \ge  2$ with smooth boundary $\partial M$, and let  $ i^*\!\!:\! \Lambda M\to \Lambda \partial M$ be the pullback via the inclusion  $i: \!\partial M \to M$.  If $\omega \in C^\infty \Lambda M$, decompose
 \begin{eqnarray*}  \omega\big|_{\partial M} =\omega^T +\omega^N \;\;\, \mbox{where}\;\; \omega^T=i^* \omega \in C^\infty \Lambda \partial M.\end{eqnarray*}  
Denote by   $\pi^T$ and $\pi^N$ the corresponding projection operators. Let
\begin{eqnarray}  \label{24.2.6-5}\mathcal{B}^a \omega = \pi^N \omega \oplus \pi^N d\omega.\end{eqnarray} 
Let $*$ be the Hodge star operator and let  $\mathcal{B}^r\omega =\mathcal{B}^a * \omega$.   As defined in \cite{BGOP}, $\omega$ is said to satisfy {\bf absolute boundary conditions} if $\mathcal{B}^a \omega =0$; $\omega$ is said to satisfy {\bf relative boundary conditions} if
$\mathcal{B}^r \omega=0$.  

We first calculate the absolute (respectively, relative) boundary conditions when $M$ is a  bounded domain $M$ in the  Euclidean space $\mathbb{R}^d$ with smooth boundary. In local coordinates,  let \begin{eqnarray} \label{24.2.4-10}\omega= u_k \,dx^k\in \Lambda^1 M.\end{eqnarray} 
Then \begin{eqnarray}\label{24.2.6-1} d\omega= \frac{\partial u_k}{\partial x^j}\,dx^j 
\wedge dx^k.\end{eqnarray}   
From  the condition $\pi^N \omega =0$, we immediately get  (here the vector fiels $\mathbf{e}_1=\frac{\partial}{\partial x_1}, \cdots, \mathbf{e}_{d-1}=\frac{\partial }{\partial x_{d-1}}$ are tangent to $\partial M$ while $\mathbf{n}=\mathbf{e}_d=\frac{\partial }{\partial x_d}$ is orthogonal to it at each point)
\begin{eqnarray} \label{24.2.5-1} u_d=0 \;\;\mbox{on}\;\; \partial M.\end{eqnarray}  And from the condition  $\pi^N d\omega=0$, we have \begin{eqnarray} 
\label{24.2.5-2} \frac{\partial u_k}{\partial x^d}=0\;\;\, \mbox{on} \;\; \partial M,  \;\;\; \mbox{for}\;\;\, k=1,\cdots, d-1.\end{eqnarray}  
Conversely, (\ref{24.2.5-1}) and (\ref{24.2.5-2}) imply $\mathcal{B}^a \omega =0$ for such a $M\subset \mathbb{R}^d$. 
Thus, the absolute boundary conditions  $\mathcal{B}^a \omega=0$ are equivalent to the conditions (\ref{24.2.5-1}) and (\ref{24.2.5-2})  for a bounded Euclidean domain $M\subset \mathbb{R}^d$. 

Next,  we find from (\ref{24.2.4-10})   that \begin{eqnarray} \label{24.2.6-3}  *\omega = (-1)^{k-1} u_k \, dx^1 \wedge \cdots \wedge dx^{k-1}\wedge \widehat{dx^k} \wedge dx^{k+1} \wedge \cdots \wedge dx^d, \end{eqnarray}
where the hat means that $dx^k$ is omitted, 
so that \begin{eqnarray*} \label{24.2.5-5}  d\,* \omega = \big(\frac{\partial u_1}{\partial x^1} +\cdots + \frac{\partial u_d}{\partial x^d}
\big) \, dx^1\wedge \cdots \wedge dx^d.\end{eqnarray*}
It follows that  \begin{eqnarray} \label{24.2.16-6}  \pi^N d*\omega =\mbox{div}\, u.\end{eqnarray}
By   $\mathcal{B}^r \omega=0$, (\ref{24.2.6-3}) and (\ref{24.2.16-6}) we get 
\begin{eqnarray} \label{24.2.5-20} u_1=\cdots =u_{d-1}=0 \;\;\, \mbox{on}\;\;\, \partial M,\end{eqnarray} 
and \begin{eqnarray} \label{24.2.5-21} \mbox{div}\; \mathbf{u}=0 \;\;\, \mbox{on}\;\;\, \partial M.  \end{eqnarray} 
Conversely,  (\ref{24.2.5-20}) and (\ref{24.2.5-21})  imply $\mathcal{B}^r \omega =0$ for $M\subset \mathbb{R}^d$. Therefore,  for a bounded domain $M$ in the Euclidean space $\mathbb{R}^d$,  the relative boundary conditions are equivalent to the conditions (\ref{24.2.5-20}) and (\ref{24.2.5-21}).

\vskip 0.25 true cm 

Now, we introduce the {\it partition function}, or  the {\it trace of the heat semigroup} for the Lam\'{e} operator $\mathscr{L}$ with mixed boundary condition $\aleph$, by
$\mathcal{Z}^\aleph (t) := \mbox{Tr}\; e^{-t \mathscr{L}^\aleph}= \sum_{k=1}^\infty e^{-t \Lambda_k^\aleph}$,  defined for $t>0$ and monotone decreasing in $t$, where $\Lambda_k^\aleph$ is the $k$-th eigenvalue of $\mathscr{L}^\aleph$ with mixed boundary conditions $\aleph$, \ $(\aleph=\{\mbox{DF}, \mbox{FD}\})$.

\vskip 0.22 true cm

\noindent{\bf Theorem 2.1 (Branson,  Gilkey,  {\O}rsted and A. Pierzchalski,  see  Theorem 4.6.4 of \cite{Gil-95}).} \ {\it Suppose  $(M,g)$ is a smooth compact  connected Riemannian manifold of dimension $d$ with smooth boundary $\partial M$. Let the Lam\'{e} coefficients $\mu$ and $\lambda$ satisfy $\mu>0$ and $\lambda+2\mu>0$, and let $0\le \Lambda_1^\aleph\le  \Lambda_2^\aleph \le \Lambda^\aleph_3\le \cdots \le \Lambda_k^\aleph \le \cdots$ be the eigenvalues of the Lam\'{e} operator $\mathscr{L}^\aleph$  with respect to the mixed  boundary conditions $\aleph$,  where  $\aleph \in
 \{\mbox{DF}, \mbox{FD}\}$. Then
\begin{eqnarray} \label{1-7} && \mathcal{Z}^{\aleph}(t) =  \sum_{k=1}^\infty e^{-t \Lambda_k^{\aleph}} =\bigg[ \frac{d-1}{(4\pi \mu t)^{d/2}}
 + \frac{1}{(4\pi (2\mu+\lambda) t)^{d/2}}\bigg] {\mbox{Vol}_d}(M) \\
&& \qquad \;\, \quad \, +\, @ \frac{1}{4} \bigg[  \frac{n-3}{(4\pi \mu t)^{(d-1)/2}}
 +  \frac{1}{(4\pi (2\mu+\lambda) t)^{(d-1)/2}}\bigg]{\mbox{Vol}_{d-1}}(\partial M)+O(t^{{1-d}/2})\quad\;\, \mbox{as}\;\; t\to 0^+,\nonumber\end{eqnarray}
where  \begin{eqnarray*} @ =\left\{\! \begin{array}{ll} - \;\;\;\;\mbox{for}\;\,\, \aleph =DF,\\
+ \;\;\;\; \mbox{for}\;\;\, \aleph=FD.\end{array}\right.  \end{eqnarray*} }

\vskip 0.2 true cm 

\begin{proof}     Since $\mathscr{L}^\aleph$ is a strongly elliptic, self-adjoint  operator for each $\aleph \in \{\mbox{DF}, \mbox{FD}\}$,  we  get that the heat semigroup $e^{-t \mathscr{L}^\aleph}$ has the following representation:  \begin{eqnarray} \label{24.1.23-1} \mathbf{G}^\aleph (t,x,y)= \sum_{k=1}^\infty  e^{-t\Lambda_k^\aleph} \mathbf{u}_k^\aleph (x) \otimes\mathbf{u}_k^{\aleph} (y), \end{eqnarray}  and its trace is 
\begin{eqnarray}  \label{24.1.23-2}  \mbox{tr}\, e^{-t  \mathscr{L}^\aleph}  = \sum_{k=1}^\infty e^{-t\Lambda_k^\aleph} = \int_{M} \mathbf{G}^\aleph (t,x,x)\, dx,\end{eqnarray} 
where $\mathbf{u}_k^\aleph$ is  eigenvector corresponding to the $k$-th eigenvalue $\Lambda_k^\aleph$.   
It follows from general principles (see e.g. \cite{Gil-95},  \cite{Gru-86} or \cite{Ta-2}) that  (\ref{24.1.23-2}) has an asymptotic expansion
\begin{eqnarray}  \label{24.1.23-3}  \mbox{tr}\, e^{-t  \mathscr{L}^\aleph}  \sim  \sum_{j=0}^\infty  b_j \,t^{(j-d)/2} \;\;\, \, \mbox{as} \ \;\,  t\to 0^+,\end{eqnarray} 
 where the coefficients $b_j$ are integrals of local expressions in the jets of the symbol of $\mathscr{L}$.  From the invariant theory in \cite{Gil-95}, it follows that  the first two terms are given by 
 \begin{eqnarray*}   b_0= \omega_0 \cdot \mbox{Vol}_d\, (M), \;\;\;\;\,  b_1= \omega_1 \cdot \mbox{Vol}_{d-1}\, (\partial M),\end{eqnarray*} 
 where $\omega_0$ and $\omega_1$ are universal constants depending only on the dimension (and $\mathscr{L}$), and 
 \begin{eqnarray*}  \omega_0= \frac{d-1}{(4\pi \mu )^{d/2}}
 + \frac{1}{(4\pi (2\mu+\lambda) )^{d/2}}.\end{eqnarray*}
 
  In order to calculate $\omega_1$ we first look at a model case, namely the $d$-dimensional Euclidean cube $M=[0,\pi ]^d\subset \mathbb{R}^d$ as in \cite{PiOr-96}.  
  Although the boundary of $M$ is not smooth, its singularities will only contribute to the high-order terms beyond the first two. 
 Clearly, the Lam\'{e} operator defined in $[0, \pi]^d$ has the following expression:
  \begin{eqnarray}\label{24.1.23-6} \!\!\!\! \mathscr{L} \mathbf{u} = \left\{ \mu\begin{pmatrix}  \Delta  & 0  & \cdots & 0\\
  0 & \Delta & \cdots & 0 \\
   \cdots & \cdots & \cdots & \cdots\\
   0 & 0 & \cdots & \Delta
    \end{pmatrix}  
    + (\lambda+\mu) \begin{pmatrix} \frac{\partial^2} {\partial x_1^2} & \frac{\partial^2 }{\partial x_1 \partial x_2} & \cdots & \frac{\partial^2}{\partial x_1 \partial x_d}\\
\frac{\partial^2} {\partial x_2 \partial x_1} & \frac{\partial^2 }{\partial x_2^2} & \cdots & \frac{\partial^2}{\partial x_2 \partial x_d}   \\
 \cdots & \cdots & \cdots & \cdots \\
 \frac{\partial^2} {\partial x_d \partial x_1} & \frac{\partial^2 }{\partial x_d \partial x_2} & \cdots & \frac{\partial^2}{\partial x_d^2}  \end{pmatrix} \right\} \begin{pmatrix} u^1\\
 u^2\\ 
 \vdots \\
 u^d \end{pmatrix}. \;\; \end{eqnarray}
For the above  cube $M=[0,\pi]^d$,  from the first  condition   of the DF boundary conditions:  
\begin{eqnarray*} \mathbf{u} - \left(g_{\alpha \beta} n^\alpha u^\beta \right) \mathbf{n}=\mathbf{u}- \langle \mathbf{n}, \mathbf{u}\rangle  \mathbf{n}=0 \;\;\;\mbox{on}\;\;\, \partial M,\end{eqnarray*}  
we immediately get that, in the local normal coordinates (recall that the vector fiels $\mathbf{e}_1=\frac{\partial}{\partial x_1}, \cdots, \mathbf{e}_{d-1}=\frac{\partial }{\partial x_{d-1}}$ are tangent to $\partial M$ while $\mathbf{n}=\mathbf{e}_d=\frac{\partial }{\partial x_d}$ is orthogonal to it at each point),  \begin{eqnarray} \label{24.1.24-11} u^1=\cdots = u^{d-1}=0\;\;\, \mbox{on}\;\;\, \partial M.\end{eqnarray}
And from the second  condition of DF boundary conditions: \begin{eqnarray*}g_{\alpha\beta} n^\alpha \left( \mathscr{T} u\right)^ \beta= \langle \mathbf{n}, \mathscr{T} \mathbf{u}\rangle =0 \;\;\;\mbox{on}\;\;\, \partial M,\end{eqnarray*} 
we obtain $\left(\mathscr{T} \mathbf{u}\right)^d=0$ on $\partial M$, 
i.e., \begin{eqnarray} \label{24.1.31-1}  \lambda \,(\mbox{div}\; \mathbf{u}) +2\mu \frac{\partial u^d}{\partial x_d}=0\;\; \mbox{on}\;\; \partial M.\end{eqnarray}   From (\ref{24.1.24-11}), we get  \begin{eqnarray*} \frac{\partial u^k}{\partial x_k}=0 \;\;\mbox{on} \;\;\partial  M, \;\;\,\;\; k=1,2,\cdots, d-1,\end{eqnarray*}
  and hence using this and (\ref{24.1.31-1})  we have    \begin{eqnarray*} (\lambda+2\mu)\frac{\partial u^d}{\partial x_d}=0\;\; \,\mbox{on}\;\;\partial M,\end{eqnarray*} 
   which implies  $\frac{\partial u^d}{\partial x_d}=0$ on $\partial M$. Combining this and (\ref{24.1.24-11}) we have \begin{eqnarray} \label{24.1.24-12} \mbox{div}\; \mathbf{u}=0\;\;\mbox{on}\;\;\, \partial M.\end{eqnarray}  Clearly, boundary conditions   (\ref{24.1.24-11}) and  (\ref{24.1.24-12})  are  just the relative  boundary  conditions (i.e., (\ref{24.2.5-20}) and (\ref{24.2.5-21}))  for the Euclidean cube $M=[0,\pi]^d$.    
 
  In \cite{BGOP},  T. Branson,  P. Gilkey,  B. {\O}rsted and A. Pierzchalski
 gave the asymptotic expansions of the heat traces $\mathcal{Z}^{(a)}(t)=\mbox{tr}\, e^{-t {P}^{(a)}}$ and $\mathcal{Z}^{(r)}(t)=\mbox{tr}\, e^{-t  {P}^{(r)}}$ for the 
  generalized Ahlfors Laplacian $P:=a\, d\delta +b\,\delta d -2b\,\mbox{Ric}$  with  absolute and relative boundary conditions for a   Riemannian manifold $(M,g)$,  respectively  (see,  (b) of Theorem 4.2  and (b) of Theorem 4.3 in \cite{BGOP}, or  (b)  and (f) of Theorem 4.6.4 in \cite{Gil-95}),  where  $a$ and $b$ are positive constants, $\mbox{Ric}$ denotes the Ricci action on $1$-forms, $d$ is the exterior differential operator, and $\delta$ the adjoint operator of $d$. 
 In  Section 8 of \cite{Liu-23}, the author of the present paper showed that the generalized Ahlfors Laplacian $P$ defined on $C^\infty T^1 M$ is just the equivalent  $1$-form representation of the  Lam\'{e} operator defined on the Riemannian manifold $(M, g)$, where $a$ and $b$ are  positive constants,  and  $\epsilon \rho$ is an  arbitrary constant multiple of the Ricci tensor.  In particular,  putting   $M=[0,\pi]^d$, we find by taking $b:=\mu, \, a:=\lambda+2\mu$  in  (b) of Theorem 4.3 (or (f) of Theorem 4.6.4 in \cite{Gil-95})   that 
 \begin{eqnarray} \label{71-7} &&  \mathcal{Z}^{(r)}(t) =  \sum_{k=1}^\infty e^{-t \Lambda_k^{(r)}} =\omega_0^{(r)} \cdot \mbox{Vol}_{d} (M)\, t^{-d/2}+ 
 \omega_1^{(r)} \cdot \mbox{Vol}_{d-1} (\partial M)\, t^{-(d-1)/2}\quad\; \; \\  [1.6mm]
 && \quad \quad +o(t^{-(d-1)/2}) \;\;\;\mbox{as}\;\;\, t\to 0^+\nonumber \end{eqnarray} 
 with \begin{eqnarray} \label{24.2.6-30} \omega_0^{(r)}\!= \!\frac{d-1}{(4\pi \mu )^{d/2}}
\! + \!\frac{1}{(4\pi (2\mu \!+\! \lambda) )^{d/2}} \;\, \mbox{and}\;\, \omega_1^{(r)}\!= \!-\frac{1}{4} \bigg(  \frac{n-3}{(4\pi \mu )^{(d-1)/2}}
 \!+  \!\frac{1}{(4\pi (2\mu \!+\!\lambda))^{(d-1)/2}}\bigg).\end{eqnarray}  
Sincel $\omega_0^{(r)}$ and $\omega_1^{(r)}$ in (\ref{24.2.6-30}) are constants depending only on the dimension $d$ and the Lam\'{e} constants $\lambda$ and $\mu$,  and since  the heat trace $\mathcal{Z}^{DF}(t) $ has asymptotic expansion (\ref{24.1.23-3}),      it follows  that  \begin{eqnarray} \label{72-7} &&  \mathcal{Z}^{DF}(t) =  \sum_{k=1}^\infty e^{-t \Lambda_k^{DF}} =\omega_0^{(r)} \cdot \mbox{Vol}_{d} (M) \,t^{-d/2}+ 
 \omega_1^{(r)} \cdot \mbox{Vol}_{d-1} (\partial M)\, t^{-(d-1)/2}\quad  \quad \; \\ 
[1.5mm] && \quad \quad\quad \quad  +o(t^{-(d-1)/2}) \;\;\;\mbox{as}\;\;\, t\to 0^+ \nonumber \end{eqnarray} 
still holds for any  compact Riemannian manifold $(M, g)$ with  boundary $\partial M$  for the above given $\omega_0^{(r)}$ and $\omega_1^{(r)}$.  The first desired result in (\ref{1-7}) is verified.

        For the Euclidean cube $[0,\pi]^d$, from the second  condition of the FD boundary conditions: \begin{eqnarray*} g_{\alpha \beta} n^\alpha u^\beta =\langle \mathbf{n}, \mathbf{u}\rangle=0 \;\;\;\mbox{on}\;\;\, \partial M, \end{eqnarray*} we get 
         \begin{eqnarray} \label{24.1.24-19} u^d=0 \;\;\;\mbox{on}\;\;\,  \partial M.\end{eqnarray}
          And from the first condition of FD boundary conditions:  $\mathscr{T}\mathbf{u} - \langle \mathbf{n}, \mathscr{T} \mathbf{u} \rangle \,\mathbf{n} =0$ on $\partial M$, we get 
         \begin{eqnarray}  \label{24.1.24-20} \frac{\partial u^1}{\partial x_d}= \cdots =\frac{\partial u^{d-1}}{\partial x_d}=0\;\; \; \mbox{on}\;\;\, \partial M.\end{eqnarray} 
         Clearly,    (\ref{24.1.24-19}) and (\ref{24.1.24-20}) exactly are the absolute  boundary conditions (i.e.,  (\ref{24.2.5-1}) and (\ref{24.2.5-2}))  for $[0,\pi]^d$.  It follows from (b) of Theorem 4.2 in \cite{BGOP} (or (b) of Theorem 4.6.4 in \cite{Gil-95}) that 
          \begin{eqnarray} \label{51-7} &&  \mathcal{Z}^{(a)}(t) =  \sum_{k=1}^\infty e^{-t \Lambda_k^{(a)}} =\omega_0^{(a)} \cdot \mbox{Vol}_{d} (M) \,t^{-d/2}+ 
 \omega_1^{(a)} \cdot \mbox{Vol}_{d-1} (\partial M)\, t^{-(d-1)/2} \\  [1.6mm]
 && \quad \quad +o(t^{-(d-1)/2}) \;\;\;\mbox{as}\;\;\, t\to 0^+\nonumber \end{eqnarray} 
 with \begin{eqnarray} \label{24.24.6-30} \omega_0^{(a)}\!= \!\frac{d-1}{(4\pi \mu )^{d/2}}
\! + \!\frac{1}{(4\pi (2\mu \!+\! \lambda) )^{d/2}} \;\, \mbox{and}\;\, \omega_1^{(a)}\!= \!\frac{1}{4} \bigg(  \frac{n-3}{(4\pi \mu )^{(d-1)/2}}
 \!+  \!\frac{1}{(4\pi (2\mu \!+\!\lambda))^{(d-1)/2}}\bigg).\end{eqnarray}  
    It is  similar to the above argument (i.e., $\omega_0^{(a)}$ and $\omega_1^{(a)}$  are the known constants depending only on $d$, $\lambda$ and $\mu$, and the heat trace $\mathcal{Z}^{FD}(t)$ has an asymptotic expansion (\ref{24.1.23-3}))  that  \begin{eqnarray} \label{72-7} &&  \mathcal{Z}^{FD}(t) =  \sum_{k=1}^\infty e^{-t \Lambda_k^{FD}} =\omega_0^{(a)} \cdot \mbox{Vol}_{d} (M) \,t^{-d/2}+ 
 \omega_1^{(a)} \cdot \mbox{Vol}_{d-1} (\partial M)\, t^{-(d-1)/2}\quad  \quad \; \\ 
[1.5mm] && \quad \quad\quad \quad  +o(t^{-(d-1)/2}) \;\;\;\mbox{as}\;\;\, t\to 0^+ \nonumber \end{eqnarray} 
still holds for any  compact Riemannian manifold $(M, g)$ with boundary $\partial M$ for the above given $\omega_0^{(a)}$ and $\omega_1^{(a)}$.  The second  desired result in (\ref{1-7}) is verified.     
    \end{proof}

 \vskip 0.48 true cm

{\it Proof of  (\ref{24.1.21-13})}.    Now, let us come back to an alternative  proof of (\ref{24.1.21-13}).   For the given Riemannian manifold $(M, g)$ with smooth boundary, since the corresponding billiards is neither dead-end nor absolutely periodic, it follows that there exist  two-term asymptotic expansions (\ref {2023.5.5-1}) for the counting functions with mixed boundary conditions. 
 Note that  the  partition function  $\mathcal{Z}^{\aleph}(t)=\sum_{k=1}^\infty e^{-t\Lambda_k^{\aleph}}$  is  just
 the Riemann-Stieltjes integral of $e^{-t\Lambda}$ with respect to the counting function ${N}^{\aleph}(\Lambda)$ 
  for each $\aleph\in \{\mbox{DF}, \mbox{FD}\}$,  i.e., 
 \begin{eqnarray} \label{2023.5.5-11} \mathcal{Z}^{\aleph}(t) = \int_{-\infty}^{+\infty} e^{-t\Lambda } d{N}^{\aleph} (\Lambda).\end{eqnarray}
It is the well-known that if the following two-term spectral asymptotic holds
 \begin{eqnarray} \label{24.1.25-20} \!\!\!\!\!\!\!\!\!\!\!\!{N}^{\aleph} (\Lambda) = a_0\,\mbox{Vol}_n (M)\, \Lambda^{d/2} + a_1^{\aleph}\, \mbox{Vol}_{d-1} (\partial M)\,  \Lambda^{(d-1)/2}+ o(\Lambda^{(d-1)/2}) \;\; \mbox{as}\;\; \Lambda\to +\infty,\end{eqnarray}
then, one immediately finds by using (\ref{2023.5.5-11}) and  (\ref{24.1.25-20}) that
\begin{eqnarray}\label{24.1.25-21} \!\!\!\!\!\!\!\!\!\! \mathcal{Z}^{\aleph}(t) = b_0\, \mbox{Vol}_n( M)\, t^{-d/2} + b_1^{\aleph}\, \mbox{Vol}_{n-1} (\partial M) \,t^{-(d-1)/2} + o(t^{-(d-1)/2})\;\; \mbox{as} \;\; t\to 0^+,\end{eqnarray}
   with   \begin{eqnarray} \label{24.1.25-25} b_0 = \Gamma \Big( 1+\frac{n}{2}\Big)\, a_0, \, \; \,\;  b_1^{\aleph} = \Gamma\Big(1+\frac{n-1}{2}\Big) \,a_1^{\aleph}. \end{eqnarray}
Hence, from (\ref{1-7}) of Theorem 2.1 and (\ref{24.1.25-21})--(\ref{24.1.25-25}) we immediately get the desired expression (\ref{24.1.21-13}). That is,  Capoferri and Mann's result (\ref{24.1.21-13}) is essentially a known result in \cite{BGOP} which  was  obtained by T. Branson, P.  Gilkey,  B. {\O}rsted and A. Pierzchalski in \cite{BGOP} (see also, Theorem 4.6.4 of \cite{Gil-95}).  $\quad \quad \quad \quad \square$

 \vskip 0.32 true cm
 
\noindent{\bf Remark 2.2.} \ {\it  In the asymptotic expansions of the Lam\'{e} operator $\mathscr{L}$, the most difficult two cases are the Dirichlet and the free boundary conditions, and the corresponding asymptotic expansions of  the heat traces have been obtained in \cite{Liu-21} (see also \cite{Liu-23}) by the author of the present paper by using double manifold method.  For the DF and FD boundary conditions, the  corresponding  asymptotic expansions become quite simpler because the required  eigen-solutions can be easily obtained by separation of variables (see,  for the same result on p.$\,$118--119,  proof of Theorem 5.1 in \cite{PiOr-96}). }

 \vskip 0.92 true cm

\section{Numerical verifications}

\vskip 0.32 true cm

  Clearly, the authors of \cite{CaMa-24} showed the known asymptotic expansions  in \cite{BGOP} (or \cite{Gil-95}) for the Lam\'{e} operator with mixed boundary conditions. However, authors of \cite{CaMa-24} could not clearly determine whether their results are correct.  They therefore  used a non-professional (ridiculous)  method, so-called the ``numerical verification''. That is,   ``predict'' a limit value $\lim_{\Lambda\to +\infty} \frac{1}{\mbox{Vol}_{d-1}(\partial M)\,\Lambda^{(d-1)/2}}\big[ {N}^{\aleph} (\Lambda)- a_0\, \mbox{Vol}_d (M) \,\Lambda^{d/2}\big]$  by only calculating the value in a finite interval $[0,3000]$ for $\Lambda$,  where $a_0$ is the one-term coefficient of
   ${N}^\aleph(\Lambda)$. 
 In fact, according to  Cauchy's limit theory, such a numerical verification can not describe any true asymptotic behaviour of $\frac{1}{\mbox{Vol}_{d-1}(\partial M)\,\Lambda^{(d-1)/2}}\big[ {N}^\aleph (\Lambda)- a \,\mbox{Vol}_d (M) \,\Lambda^{d/2}\big]$  as $\Lambda\to +\infty$. In principle, numerical experimental verification should consider the case for sufficiently large $\Lambda$.  More precisely,  in order to verify $$\lim_{\Lambda\to +\infty} \frac{1}{\mbox{Vol}_d\,(\partial M)\,\Lambda^{(d-1)/2} }\left[ {N} (\Lambda)- a_0 \mbox{Vol}_d (M) \,\Lambda^{d/2} \right] = b'$$ for some constant $b'$, by Cauchy's limit definition, one should prove that for every number $\epsilon>0$ there is a number $M>0$ which depends on $\epsilon$ such that if $\Lambda>M$ then $$\bigg|\frac{1}{\Lambda^{(d-1)/2}}\left[ {N} (\Lambda)- a_0 \mbox{Vol}_d (M) \,\Lambda^{d/2} \right]  - b'\bigg|<\epsilon.$$
   However, the authors  of {\rm\cite{CaMa-24}} did not use such a standard method. Instead, the  ``verify''  their result by using (non-professional)   ``numerically'' verifying for $0\le \Lambda\le 3000$  for the unit disk and flat cylinders. {\bf Such a so-called ``numerical verification'' method is not believable at all. The more serious problem is that  their eigenvalues calculations in  \cite{CaMa-24}  were based on some erroneous formulas for a unit disk and flat cylinders} (see below remark).

\vskip 0.26 true cm

 On p.18 of \cite{CaMa-24}, the authors wrote: \textcolor{blue}{\it  ``Let $M\subset \mathbb{R}^2$  be the unit disk and let us work in standard polar coordinates $(r,\theta)$. Following Chapter XIII of \cite{MoFe-53}  (see also \cite{LeMoSe-21}), we introduce a fictitious third coordinate z orthogonal to the disk and seek solutions in the form
\begin{eqnarray}\label{24.1.25-31} \mathbf{u}(r,\theta) = \mbox{grad}\;\psi_1(r,\theta) +\mbox{curl}\,  (\psi_2(r,\theta)\,\mathbf{\hat{z}}) \qquad \;\; \;\qquad \quad \end{eqnarray}
 where $\mathbf{\hat{z}}$ is the unit vector in the direction of  $z$ and $\psi_j$, $j=1,2$,  are auxiliary scalar potentials. Substituting (\ref{24.1.25-31})  into  $\mathscr{L} \mathbf{u} =\Lambda\mathbf{u}$ one obtains that the scalar potentials must satisfy the Helmholtz equations
\begin{eqnarray} \label{24.1.25-32}   & - \Delta \psi_j = \omega_{j, \Lambda}  \psi_j, \;\;\;\;\; j=1,2,  \\
 \label{24.1.25-33}&  \omega_{1, \Lambda} :=\frac{\Lambda}{\lambda+2\mu},\;\; \;\;\,  \omega_{2, \Lambda}= \frac{\Lambda}{\mu}. \end{eqnarray} 
But now the general solution to (\ref{24.1.25-32}) regular at $r = 0$ reads
\begin{eqnarray} \label{24.1.25-34}  \psi_j(r, \phi)= c_{j,0}J_0(\sqrt{\omega_{j,\Lambda}}\,r)+\sum_{k=1}^\infty J_k(\sqrt{\omega_{j,\Lambda} }\, r)\big(c_{j,k,+} e^{i k\phi} + c_{j,k,-} e^{-ik\phi}\big), \;\;\quad  \end{eqnarray}
where the $J_k$'s are Bessel functions of the first kind.''
}

\vskip 0.25 true cm

The above equations (\ref{24.1.25-32}) are incorrect. The correct statement should be (cf. Theorem 3.1 below)
$$ -  \Delta (\mbox{grad}\,\psi_1) = \omega_{1, \Lambda} (\mbox{grad}\, \psi_1),\;\;\; \;\;  - \Delta \big(\mbox{curl}\, (\mathbf{z} \psi_2)\big)= \omega_{2, \Lambda} \big(\mbox{curl}\, (\mathbf{z} \psi_2)\big)$$
instead of (\ref{24.1.25-32}).
That is, \begin{eqnarray*} \mbox{grad}\,\big( \Delta \psi_1 + \omega_{1, \Lambda} \psi_1\big)=0,\;\;\; \;\;  \mbox{curl}\,\big( \Delta  (\mathbf{z} \psi_2)+ \omega_{2, \Lambda} \, (\mathbf{z} \psi_2)\big)=0\end{eqnarray*}
because of  \begin{eqnarray*} \Delta \circ \mbox{grad} =\mbox{grad}\circ \Delta\;\;\,\mbox{and}\;\;  \Delta \circ \mbox{curl}= \mbox{curl} \circ \Delta. \end{eqnarray*} 
Or equivalently,
\begin{eqnarray} \label{2023.5.14-1} &\Delta \psi_1 + \omega_{1, \Lambda}  \psi_1=C \;\;\mbox{for any constant}\;\;C\in \mathbb{R}^1,\\
  \label{2023.5.14-2} &\Delta ( \mathbf{z} \psi_2)+ \omega_{2, \Lambda} \, (\mathbf{z} \psi_2)=\mathbf{f}, \;\,\mbox{for any} \;\; \mathbf{f}\;\; \mbox{with}\;\; \mbox{curl}\,\mathbf{f}=0.\end{eqnarray}
 In \cite{CaMa-24}, only special  $C\equiv 0$ and $\mathbf{f}\equiv 0$ are chosen. 
 Obviously, for any $C\ne 0$ or $\mathbf{f}\ne 0$ with $\mbox{curl}\ \mathbf{f}=\mathbf{0}$,  the equations (\ref{2023.5.14-1}) and (\ref{2023.5.14-2})  are non-homogenous equations, their solutions  have different forms except for  those in (\ref{24.1.25-34}).
 Thus, a large number of solutions have not been considered in  \cite{CaMa-24}, and a large number of elastic eigenvalues have been lost.
{{\bf Because  the above  (\ref{24.1.25-32})  (by \cite{CaMa-24}) \ is incorrect,  all calculations for elastic eigenvalues in the unit disk (in particular, (\ref{2023.5.14-1}) and (\ref{2023.5.14-2}))  in  \cite{CaMa-24} are wrong. Of course, Fig.$\,$1 and Fig.$\,$2 (on p.$\,$ 18-19 in \cite{CaMa-24})  are also wrong.}}

\vskip  0.25 true cm
{{\bf This mistake stems from the supplementary of an earlier erroneous papers \cite{LeMoSe-21},   and once more the same mistake occurred in
 \cite{CaFrLeVa-23}  and \cite{CaMa-24}.
}}

\vskip 0.45 true cm
To help with a good understanding to the above discussions, here we copy 2.5.Theorem on p.$\,$123--124 of \cite{KGBB}:

\vskip  0.2 true cm
    \noindent{\bf Theorem 3.1.} \ {\it  Let $D$ be a domain in $\mathbb{R}^3$. The solution $\mathbf{u} = (u^1,u^2,u^3)\in [C^2(D)\cap C^1(\bar D)]^3$ of equation \begin{eqnarray} \label{2023.5.6-3}  \mu \Delta\mathbf{u}+(\lambda+\mu)\, \mbox{grad} \; \mbox{div}\, \mathbf{u}+\Lambda \mathbf{u}=0 \;\, &\mbox{in}\;\, D,\end{eqnarray} is
represented as the sum
\begin{eqnarray} \label{2023.5.6-6} \mathbf{u}=\mathbf{u}^{(p)} +\mathbf{u}^{(s)},\end{eqnarray}
where $\mathbf{u}^{(p)}$  and $\mathbf{u}^{(s)}$  are the regular vectors, satisfying the conditions
\begin{eqnarray}\label{2023.5.6-7} (\Delta +\omega_{1,\Lambda})\mathbf{u}^{(p)} =0, \;\; \mbox{curl}\; \mathbf{u}^{(p)}=0,\\
 (\Delta +\omega_{2,\Lambda})\mathbf{u}^{(s)} =0, \;\; \mbox{div}\; \mathbf{u}^{(s)}=0, \;\;\end{eqnarray}
where
\begin{eqnarray} \omega_{1,\Lambda} := \frac{\Lambda}{ \lambda+2\mu}, \;\;\, \omega_{2,\Lambda} := \frac{\Lambda}{ \mu} \end{eqnarray}
}

\vskip 0.35 true cm
\noindent{\bf Remark 3.2.} \  {\it The similar mistake also occurred in the example of flat  cylinders in \cite{CaMa-24}. Of course, all  discussions and all figures from p.$\,$19 to p.$\,$25 in  \cite{CaMa-24} are all wrong,  we  refer the reader to \cite{Liu-23} for the detailed remarks.}

\vskip 0.35 true cm
\noindent{\bf Remark 3.3.} \  {\it  It is  quite strange:  why does correctness of a correct (proved) known conclusion in \cite{CaMa-24} need to be verified by a flawed ``numerical verification''?  More worst, these calculation formulae of  numerical verifications  in \cite{CaMa-24} are all wrong.}

\vskip 0.45 true cm

\vskip 1.26 true cm


\begin{thebibliography}{99}

\vskip 0.39 true cm 





\bibitem[BGOP]{BGOP}   T. P. Branson, P. B. Gilkey, B. {\O}rsted,  A. Pierzchalski, {\it Heat equation asymptotics of a
generalized Ahlfors Laplacian on a manifold with boundary},  Operator Calculus and Spectral Theory,  M. Demuth, B. Gramsch and B. Gramsch, and B.Schulze,  eds., p$\,$. 1--13,  Birkh$\mathrm{\ddot a}$user, Boston, 1992.


\bibitem[CaMa-24]{CaMa-24}  M. Capoferri and I.  Mann,  \textsl{Spectral asymptotics for linear elasticity:
the case of mixed boundary conditions},  arXiv:2311.18652.


\bibitem[CaFrLeVa-23]{CaFrLeVa-23}
M. Capoferri,  L. Friedlander, M. Levitin,  D. Vassiliev, {\it Two-Term Spectral Asymptotics in Linear Elasticity},
The Journal of Geometric Analysis (2023) 33:242, https://doi.org/10.1007/s12220-023-01269-y.



 \bibitem[Gil-95]{Gil-95} P. Gilkey, {\it Invariance Theory, the Heat Equation and the Atiyah--Singer Index Theorem}, CRC Press, Boca Raton, 1995.

\bibitem[Gil-75]{Gil-75}  P.  Gilkey, \textsl{The spectral geometry of a Riemannian manifold}, J. Differential Geometry, 10(1975), 601--618.


   \bibitem[Gru-86]{Gru-86} G. Grubb, {\it Functional Calculus of Pseudo-Differential Boundary Problems}, Birkh\"{a}user, Boston, 1986.
    
    
    \bibitem[Ho4]{Ho4}  L. H\"{o}rmander, {\it The Analysis of Partial Differential Operators IV},  Springer-Verlag, Berlin Heidelberg New York, 1985.

    


\bibitem[KGBB]{KGBB}  V. Kupradze, T. Gegelia, M. Basheleishvili, T. Burchuladze, \textsl{Three-dimensional problems of the mathematical theory
of elasticity and thermoelasticity}, North-Holland Publishing Company, Amsterdam, New York, Oxford 1979.


  
  \bibitem[LeMoSe-21]{LeMoSe-21}  M. Levitin, P. Monk,  V. Selgas, {\it Impedance eigenvalues in linear elasticity}, SIAM J. Appl. Math., 81:6(2021), 2433-2456.


\bibitem[LiQin-13]{LiQin-13}  T. Li and T. Qin,  {\it Physics and partial differential equations}, Vol.1, Higher Education Press, Beijing, 2013, Translated by Y. Li.


\bibitem[Liu-21]{Liu-21} G. Q. Liu, {\it Geometric Invariants of Spectrum of the Navier--Lam\'{e} Operator}, J. Geom. Anal. 31 (2021), 10164--10193.

\bibitem[Liu-23]{Liu-23} G. Q. Liu, {\it Remark on paper ``Two-term spectral asymptotic in linear elasticity''}, arXiv: 2305.14134v6.

 \bibitem[Liu-19]{Liu-19} G. Q. Liu, {\it Determination of isometric real-analytic metric and spectral invariants for elastic Dirichlet-to-Neumann map on Riemannian manifolds}, arXiv:1908.05096.


\bibitem[Liu-22]{Liu-22} G. Q. Liu, {\it  Two-term spectral asymptotics in linear elasticity on a Riemannian manifold},  arXiv:2208.02679v6.   

    \bibitem[Liu-2309]{Liu-2309} G. Q. Liu,  \textsl{On an algorithm for two-term spectral asymptotic formulas}, arXiv:2309.07475v3.


\bibitem[MaHu]{MaHu} J. E. Marsden and T. R. Hughes, {\it Mathematical Foundations of elasticity,  Dover Publications}, Inc., New York, 1983.


\bibitem[MoFe-53]{MoFe-53} P. M. Morse and H. Feshbach, {\it 
 Methods of theoretical physics} , Vol. 2, McGraw-Hill, N. Y., 1953.

\bibitem[PiOr-96] {PiOr-96}  A. Pierzchalski and  B. {\O}rsted, \textsl{The Ahlfors Laplacian on a Riemannian manifold with boundary},
 Michigan Math. J., 4(1996). 

    
 \bibitem[SaVa-97]{SaVa-97} Yu. Safarov and D. Vassiliev, {\it The asymptotic distribution of eigenvalues of partial differential
operators}, Amer. Math. Soc., Providence, RI, 1997. DOI: 10.1090/mmono/155.



\bibitem[Ta-2]{Ta-2} M. E. Taylor, {\it Partial Differential
Equations II}, \ 2nd Edition, Appl. Math. Sci., vol. 116, Springer Science+Business Media, LLC 1996, 2011.


\bibitem[Ta-3]{Ta-3} M. E. Taylor, \textsl{Partial Differential
Equations III}, Second Edition, Appl. Math. Sci., vol. 117, Springer-Verlag, New York, 2011.

 
    
\bibitem[Va-84]{Va-84}   G. G. Vasil'ev, {\it  Two-term asymptotics of the spectrum of a boundary value problem under an interior reflection of general form}. Funkts. Anal. Pril. 18(4), 1-13 (1984) (Russian, full text available at Math-Net.ru); English translation in Funct. Anal. Appl. 18, 267-277 (1984). https://doi.org/10.1007/BF01083689.

\bibitem[Va-86]{Va-86} D. G. Vasil'ev, {\it Two-term asymptotic behavior of the spectrum of a boundary value problem in the case of a piecewise smooth boundary}, Dokl. Akad. Nauk SSSR 286(5), 1043--1046 (1986) (Russian, full text available at Math-Net.ru); English translation in Soviet Math. Dokl. 33(1), 227--230 (1986), full text available at the author’s website https://www.ucl.ac.uk/~ucahdva/publicat/vassiliev86.pdf.




\end{thebibliography}
\end{document}